# Applying Quantum Computing to Solve Multicommodity Network Flow Problem

Niu Chence*
School of Civil and Environmental Engineering, UNSW Sydney
chence.niu@unsw.edu.au

Purvi Rastogi
Delhivery Limited, India
purvi.rastogi@delhivery.com

Jaikishan Soman
Delhivery Limited, India
s.jaikishan@delhivery.com

Kausik Tamuli
Delhivery Limited, India
kausik.tamuli@delhivery.com

Vinayak V. Dixit
School of Civil and Environmental Engineering, UNSW Sydney
v.dixit@unsw.edu.au

*Corresponding author

**Abstract**

In this paper, the multicommodity network flow (MCNF) problem is formulated as a mixed integer programing model which is known as NP-hard, aiming to optimize the vehicle routing and minimize the total travel cost. We explore the potential of quantum computing, specifically quantum annealing, by comparing its performance in terms of solution quality and efficiency against the traditional method. Our findings indicate that quantum annealing holds significant promise for enhancing computation in large-scale transportation logistics problems.

1. Introduction

In the context of urban express parcel transportation, the challenge lies in optimizing a complex distribution network that encompasses multiple pick-up, processing, and delivery sites. This network is integral to multicommodity network flow (MCNF) Problems, where each site and the routes connecting them must be efficiently organized [1]. Key influencing factors include service timelines, load volumes, and vehicle operational costs, all of which significantly impact daily operational expenses—a substantial component of total costs.

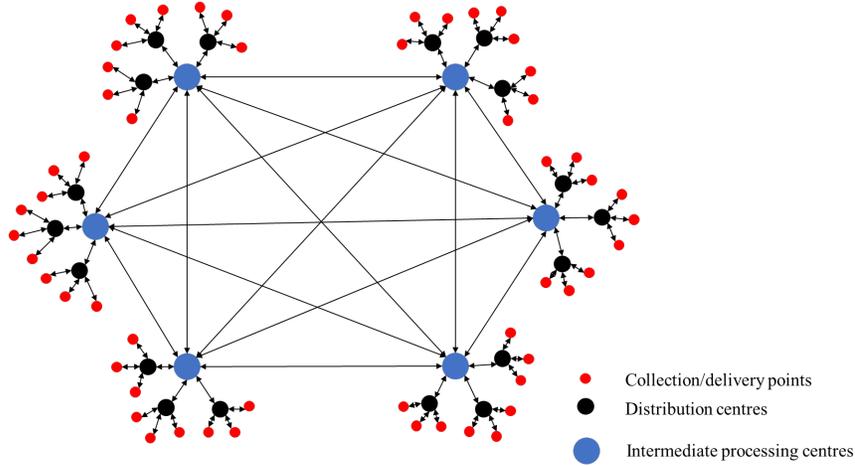

Figure 1 Topology of the distribution network

Figure 1 depicts the simplified structure of the distribution network for a logistics service provider. The network's initial phase encompasses pickup centers (PCs) and distribution centres (DC) which are often co-located. Most pickup load is brought to these locations and then transferred to the corresponding intermediate processing centers (IPCs). These origin IPCs serve as aggregation points for all shipments gathered from nearby DCs/PCs. Based on the destination, shipments from the origin IPC are routed to the destination IPC, either directly or via multiple IPCs. The diagram illustrates the links between IPCs, known as the 'primary line haul'. Shipments are then conveyed from the destination IPC to the DC via the secondary linehaul network.

For each shipment, the originating and destination IPCs are identifiable, allowing for the calculation of load transfer between any two IPCs. The goal of MCNF is to identify least cost paths that efficiently convey loads from each origin to destination IPC, in compliance with network constraints. Associated with each origin and destination IPC, we also have the total turnaround time (TAT) that is available for the transfer of the loads. In this paper, we propose a mixed-integer programming model for the MCNF problem. By using quantum computers, more specifically quantum annealing, these problems are solved by the hybrid constrained quadratic model (CQM) solver provided by D-Wave.

In the next sections, we briefly introduce quantum computing as well as MCNF problems. Next,

the mathematical statement of the problem is proposed. Lastly, the quantum computing capability is illustrated by comparing to the conventional CPLEX solver.

## 2. Literature review

This section provides an overview of the present statues of quantum computing, with a focus on quantum annealing. Additionally, it highlights key advancements in MCNF problems.

*2.1 Quantum computing*

Over the past three decades, quantum computing has made remarkable advancements. Its superiority over classical computing lies in the quantum mechanical properties of qubits, specifically superposition and entanglement. Unlike classical bits that exist in a binary state of either 0 or 1, a quantum bit or qubit can simultaneously embody both states due to superposition. Consequently, two qubits can represent four states: (0,0), (1,0), (0,1), and (1,1). This attribute allows for exponential growth in potential energy states with increasing qubits, granting quantum computing a significant edge over classical methods. For instance, a 20-qubit quantum computer can concurrently represent approximately 1 million (2^20) different states. A notable achievement in this field is the demonstration of 'Quantum Supremacy', where a 53-qubit Sycamore processor solved a problem in 200 seconds that would take a classical supercomputer approximately 10,000 years [2].

There are two leading paradigms of quantum computing, one is gate-based quantum computing (also known as circuit-based quantum computing) which is based on the quantum logic gates and the other one is quantum annealing which is based on the quantum tunnelling [3]. The gate-based quantum computing is similar to the classic logic gates but with quantum properties such as superposition and entanglement. Since the qubits is very sensitive to the environmental disturbances, the gate-based quantum computing faces several challenges such as quantum decoherence, error correction and qubits controlling [4]. On the other hand, the quantum annealing provide a relatively easier way to scale for specific applications such as optimization problems and probabilistic sampling problems [5]. The process of quantum annealing starts at initializing the qubits in a superposition state and gradually evolves the system to minimize a specific energy function which is defined as a following time-dependent Ising Hamiltonian:

$$H_{QA}/h = \underbrace{-A(t/t_f)\sum_i \sigma_i^X}_{Initial\ Hamiltonian} + \underbrace{B(t/t_f)\left(\sum_i h_i\sigma_i^Z + \sum_{i>j} J_{ij}\sigma_i^Z\sigma_j^Z\right)}_{Final\ Hamiltonian} \quad (1)$$

where $t_f$ is the annealing time, $\sigma_i^X$ and $\sigma_i^Z$ are the Pauli matrices acting on qubit $i$, and $h_i$ and $J_{ij}$ are the qubit biases and coupling strengths, respectively. The quantum annealing begins with the ground state of the initial Hamiltonian, where all qubits are in a superposition of $|0\rangle$ and $|1\rangle$. Following a preset annealing schedule given by the time dependent functions $A(t/t_f)$ and $B(t/t_f)$, the Hamiltonian of the system slowly changes from the initial to the final Hamiltonian state, which encodes the solution of the given optimization problem. The Quantum Annealer solves

Ising minimization problems, which are isomorphic to a Quadratic Unconstrained Binary Optimization (QUBO) Problem that are NP-Hard problems of the form:

$$\text{Obj} \coloneqq x^T Q x \qquad (2)$$

where x is a vector of N binary variables and Q is a NxN matrix representing the coefficients of the quadratic terms. The diagonal terms of Q are mapped to $h_i$ and the cross terms are mapped to $J_{ij}$ in the final Hamiltonian.

Due to the limitation of the number of qubits, hybrid-solvers are developed by D-Wave. It is a combination of classical solvers and QPU. Specifically, the problem is computed by a number of heuristic solvers. Then the quantum models map the mathematical formulation onto QPU which guides the heuristic solvers to improve the solution quality [6]. Currently, there are more studies investigating the computing capability of the quantum computing in different fields such as transportation [7–9], engineering design [10], finance [11].

*2.2 Multicommodity network flow problem*

Ford and Fulkerson [12] and Hu [13] were the first to model the MCNF problems, incorporating capacity constraints on arcs or links.. More recently, Salimifard and Bigharaz [14] presents the comprehensive review of MCNF classifications and related solution methods. The significance of MCNF is underscored by its diverse applications across various domains. It plays a critical role in transportation logistics [15–18], economic sector [17,19,20] and communication networks [21–23]. From an algorithm standpoint, approaches to solving MCNF problems are broadly classified into three categories: exact (such as column generation), approximation (such as Lagrangian relaxation) and heuristics/metaheuristics (such as genetic algorithms and simulated annealing). For example, Sarubbi et al. [24] designed a cut-and-branch algorithm based on Benders decomposition to solve the multicommodity traveling salesman problem. The numerical results showed that the designed algorithm was faster than stand-alone CPLEX. A branch-and-price-and-cut algorithm was proposed by Barnhart et al. [25] to solve the origin-destination integer MCNF problems. Retvdri et al. [26] employed the Lagrangian-relaxation technique to get an initial feasible solution for the minimum cost MCNF problem, with a focus on open shortest path first traffic engineering. Moshref-Javadi and Lee [27] proposed a heuristic method by adapting and synthesizing simulated annealing and variable neighbourhood search for MCNF vehicle routing problem.

Quantum computing is still in the Noisy Intermediate Scale Quantum (NISO) era, quantum systems are particularly vulnerable to environmental noise and errors [28,29]. Despite these challenges, given the unique properties of quantum computing, such as superposition and entanglement, quantum computing still gives the promising solutions for solving complex problems. Therefore, it is still worth investigating the computational potential of quantum computing for the large scale MCNF problem.

3. **Problem formulation**

In this section, an optimization model for the efficient transfer of loads among IPCs is introduced. In the graph $G(O,A)$, $O$ represents a set of IPC locations, while $A$ denotes the arcs connecting these locations. Let $O = 1, 2, \ldots, n$ be the set of $n$ IPCs. The graph includes $n \times (n-1)$ origin-destination (OD) pairs, collectively defined as set $K = k_1, k_2, \ldots, k_{n(n-1)}$. For each OD pair $k$ in $K$, there is an associated non-negative load $L_k$ that needs to be transferred from the corresponding origin $o_k \in O$ and destination $d_k \in O$.

The objective is to determine the required number of vehicles on each arc, to efficiently transport all loads across the network while minimizing total transportation costs. Due to the limitation of the number of qubits and the integrated control errors in quantum annealing, only one type of vehicle is used in this study. We assume that there is no upper limit on the number of vehicles available. The mixed-integer programming model is given by equations (3)-(5). The notations are defined in Table 1.

$$\min \sum_{(i,j) \in A} d_{i,j} * CV * N_{ij} \tag{3}$$

$$\sum_{j \in O} x_{kij} * L_k - \sum_{j \in O} x_{kji} * L_k = b_{ik} \forall i \in O, k \in K \tag{4}$$

$$\sum_{k \in K} x_{kij} * L_k \leq \sum_{v \in V} N_{ij} * W \ \forall \ i \in O, j \in O : i \neq j \tag{5}$$

Table 1 Mathematical notation

| | |
|---|---|
| $W$ | vehicle capacity |
| $CV$ | vehicle variable cost per unit distance |
| $dd_{ij}$ | distance in km between IPC $i$ and $j$ |
| $L_k$ | load for OD pair $k$ |
| $b_{ki} =$ | $\begin{cases} L_k, if\ i\ is\ the\ orgin\ IPC \\ -L_k, if\ i\ is\ the\ destination\ IPC \\ \ \ \ \ 0, otherwise \end{cases}$ |
| $x_{kij}$ | binary variable of load k transferred between arc $(i,j) \in A$ |
| $N_{ij}$ | the number of vehicles connecting arc $(i,j) \in A$. |

The objective function is defined as Equation (3) to minimize the total travel cost to transport load on arcs $(i,j) \in A$. Equation (4) ensures that the flow of a given load $k \in K$ is continuous from the source to the destination IPC. In hop cites $b_{ki} = 0$, so that flow coming into the city should flow out completely. Equation (5) ensures that the load transfer in an arc is less than the combined weight capacity of the vehicles assigned to the arc $(i,j) \in A$.

4. **Numerical experiments**

In this section, numerical experiments are conducted to compare the computational performance between CQM solver and single-thread CPLEX solver on a Win64 laptop with Intel (R) Core (TM) i7-1365U CPU and 32 GB RAM.

To reduce the size of graph $G$, we applied constraints to the arc set $A$. For each origin-destination

(OD) pair, a set of potential IPC locations was identified, which can be directly connected from the origin IPC. Utilizing this set of restricted connections, defined by a restriction matrix, we cataloged all feasible paths for each OD pair. Subsequently, paths that do not comply with the turnaround time (TAT) conditions were excluded. This process led to the creation of a specific subset of arcs for each OD pair, derived from these enumerated paths. Therefore, instead of considering every IPC-to-IPC connection in the network, we focused on a tailored subset of arcs for each OD pair. The implementation of the restriction matrix significantly reduced the size of graph $G$.

## 5. Results

As shown in Table 2 and Figure 2, in evaluations across varying number of variables, CPLEX outperforms hybrid-CQM solver when variables are less than 40,000. When the number of variables is between 50,000 to 200,000, hybrid-CQM solver consistently outperforms CPLEX in terms of time efficiency and performance. Specifically, when the number of variables is 20455, both methods take 5 seconds, but CPLEX's objective value is 2.14E+06 compared to CQM's 6.93E+06. As the number of variables increase (56819 to 191886), hybrid-CQM solver consistently achieves better performance values in shorter durations. For instance, at 95108 variables, hybrid-CQM solver reaches 1.62E+07 in 50 seconds, while CPLEX takes ten times longer (500 seconds) for a higher value of 3.70E+07. Beyond 150,000 variables, hybrid-CQM solver's time efficiency remains better than CPLEX, but its performance advantage lessens.

Table 2 Computational results of different size of benchmark networks

| No. of variables | No. of constraints | Hybrid-CQM solver | Single-thread CPLEX | Hybrid-CQM solver running time (s) | Single-thread CPLEX running time (s) | MIPGap* |
|---|---|---|---|---|---|---|
| 20455 | 8992 | 6.93E+06 | 2.14E+06 | 5 | 5 | 32.80% |
| 40936 | 17070 | 1.05E+07 | 2.83E+06 | 12 | 12 | 35.78% |
| 56819 | 24198 | 1.35E+07 | 3.21E+07 | 19 | 150 | 89.91% |
| 76680 | 32653 | 1.52E+07 | 3.55E+07 | 30 | 170 | 90.03% |
| 95108 | 40644 | 1.62E+07 | 3.70E+07 | 50 | 500 | 90.08% |
| 113056 | 48301 | 1.77E+07 | 3.88E+07 | 55 | 700 | 89.74% |
| 133395 | 56422 | 1.96E+07 | 4.56E+07 | 73 | 500 | 89.48% |
| 152558 | 152558 | 2.19E+07 | 2.32E+07 | 100 | 450 | 71.85% |
| 191886 | 80538 | 2.44E+07 | 2.57E+07 | 135 | 1800 | 69.76% |

MIPGap means the relative difference between the best integer solution (upper bound) and the best possible solution (lower bound) in CPLEX [30].

The computational process of CPLEX and hybrid-CQM are examined. It is found that, over a certain period, there is no significant improvement in the solution quality for both methods. For example, with CPLEX, when the number of variables is 113,056, a 94.76% gap between the best solution and bound initially exists at 0.00 seconds. As illustrated in Table 2, the MIPGap narrows only marginally to 89.74% (3.88E+07) after 700 seconds, suggesting potential stagnation or difficulties in the optimization process. It should be noted that when the computation time is approximately 720 seconds, the MIPGap is observed to be 17.37%. The challenges are observed with the hybrid-CQM solver. After 55 seconds, the output is recorded at 1.77E+07, and after 500 seconds, it improves slightly to 1.70E+07, an improvement of 3.95%. This indicates that under computational time constraints, the hybrid-CQM solver is both efficient and reliable. However, in the absence of computational time restriction, CPLEX tends to provide better solutions.

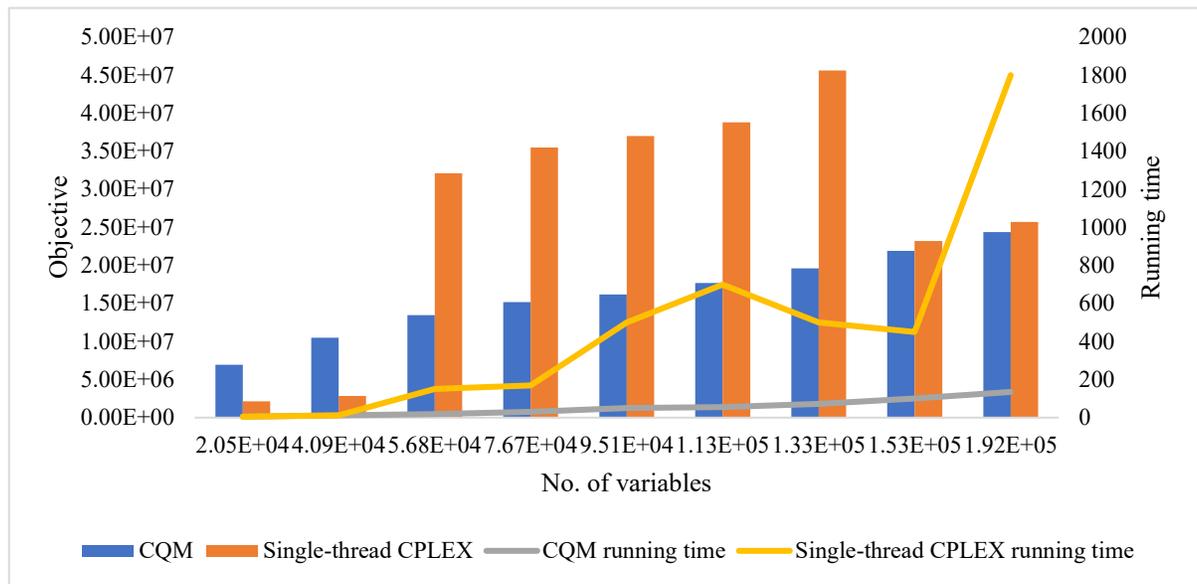

Figure 2 Computational experience of benchmark networks

## 6. Conclusion and discussion

In this paper, we present the first application of the hybrid-CQM to address MCNF problem. We proposed a mixed-integer programming model with the aim of minimizing the total travel cost. The computational performance of the hybrid-CQM solver is evaluated with the CPLEX solver across different scales of benchmark networks. The results demonstrate that when the problem size is small, the CPLEX solver shows better performance on both solution quality and computational efficiency. With the increase of number of variables, the hybrid-CQM solver shows the advantages in both aspects. Notably, the computational time of CPLEX is approximately 4.5 to 13 times longer than that of hybrid-CQM solver and the objective values ratios range approximately 1.05 to 2.39. It is worth noting that when the number of variables is larger than around 150000, the solution quality of the hybrid-CQM solver aligns closely with that of CPLEX. The possible reason could be that as the size of problem grows, there is an increased percentage of conventional heuristic algorithm in the hybrid-CQM solver.

Due to the limitation of the number of qubits and the integrated control errors, the quantum annealing is still in the early stage. Nevertheless, a feasible and practical advantage of quantum annealing is observed in solving the large-scale MCNF problem defined in the paper. With more qubits and the development of error correction techniques, quantum computing holds the significant promising computational capability to solve the large-scale transportation problems.

**Author contributions**
The authors declare no competing interests.

**Data availability statement**
Due to the confidential nature of the data used in this study, the datasets are not publicly available.

# References

[1] Wang I-L. Multicommodity network flows: A survey, Part I: Applications and Formulations. Int J Oper Res 2018;15:145–53.

[2] Arute F, Arya K, Babbush R, Bacon D, Bardin JC, Barends R, et al. Quantum supremacy using a programmable superconducting processor. Nature 2019;574:505–10.

[3] Neumann NM, de Heer PB, Phillipson F. Quantum reinforcement learning: Comparing quantum annealing and gate-based quantum computing with classical deep reinforcement learning. Quantum Inf Process 2023;22:125.

[4] Leymann F, Barzen J. The bitter truth about gate-based quantum algorithms in the NISQ era. Quantum Sci Technol 2020;5:044007.

[5] Yulianti LP, Surendro K. Implementation of quantum annealing: A systematic review. IEEE Access 2022.

[6] D-wave. Hybrid Solver for Constrained Quadratic Models. 2021.

[7] Dixit V, Jian S. Quantum Fourier transform to estimate drive cycles. Sci Rep 2022;12:1–10.

[8] Dixit VV, Niu C. Quantum computing for transport network design problems. Sci Rep 2023;13:12267. https://doi.org/10.1038/s41598-023-38787-2.

[9] Dixit VV, Niu C, Rey D, Waller ST, Levin MW. Quantum computing to solve scenario-based stochastic time-dependent shortest path routing. Transp Lett 2023:1–11. https://doi.org/10.1080/19427867.2023.2238461.

[10] Chen Z-Y, Xue C, Chen S-M, Lu B-H, Wu Y-C, Ding J-C, et al. Quantum approach to accelerate finite volume method on steady computational fluid dynamics problems. Quantum Inf Process 2022;21:137.

[11] Bova F, Goldfarb A, Melko RG. Commercial applications of quantum computing. EPJ Quantum Technol 2021;8:2. https://doi.org/10.1140/epjqt/s40507-021-00091-1.

[12] Ford LR, Fulkerson DR. A Suggested Computation for Maximal Multi-Commodity Network Flows. Manag Sci 1958;5:97–101. https://doi.org/10.1287/mnsc.5.1.97.

[13] Hu TC. Multi-Commodity Network Flows. Oper Res 1963;11:344–60. https://doi.org/10.1287/opre.11.3.344.

[14] Salimifard K, Bigharaz S. The multicommodity network flow problem: state of the art classification, applications, and solution methods. Oper Res 2022;22:1–47. https://doi.org/10.1007/s12351-020-00564-8.

[15] Dupas R, Taniguchi E, Deschamps J-C, Qureshi AG. A Multi-commodity Network Flow Model for Sustainable Performance Evaluation in City Logistics: Application to the Distribution of Multi-tenant Buildings in Tokyo. Sustainability 2020;12:2180. https://doi.org/10.3390/su12062180.

[16] Erera AL, Morales JC, Savelsbergh M. Global intermodal tank container management for the chemical industry. Transp Res Part E Logist Transp Rev 2005;41:551–66. https://doi.org/10.1016/j.tre.2005.06.004.

[17] Kuiteing AK, Marcotte P, Savard G. Pricing and revenue maximization over a multicommodity transportation network: the nonlinear demand case. Comput Optim Appl

2018;71:641–71. https://doi.org/10.1007/s10589-018-0032-0.

[18] Hernández-Pérez H, Salazar-González J-J. The multi-commodity one-to-one pickup-and-delivery traveling salesman problem. Eur J Oper Res 2009;196:987–95. https://doi.org/10.1016/j.ejor.2008.05.009.

[19] Manfren M. Multi-commodity network flow models for dynamic energy management – Mathematical formulation. Energy Procedia 2012;14:1380–5. https://doi.org/10.1016/j.egypro.2011.12.1105.

[20] Singh I, Ahn CY. A dynamic multi-commodity model of the agricultural sector: A regional application in Brazil. Eur Econ Rev 1978;11:155–79.

[21] Wadie CS, Ashour ME. Multi-commodity flow, multiple paths load balanced routing in Ad-Hoc Networks 2013.

[22] Lozovanu D, Fonoberova M. Optimal Dynamic Multicommodity Flows in Networks. Electron Notes Discrete Math 2006;25:93–100. https://doi.org/10.1016/j.endm.2006.06.087.

[23] Padmanabh K, Roy R. Multicommodity flow based maximum lifetime routing in wireless sensor network. 12th Int. Conf. Parallel Distrib. Syst. - ICPADS06, Minneapolis, MN, USA: IEEE; 2006, p. 8 pp. https://doi.org/10.1109/ICPADS.2006.67.

[24] Sarubbi J, Miranda G, Luna HP, Mateus G. A Cut-and-Branch algorithm for the Multicommodity Traveling Salesman Problem. 2008 IEEE Int. Conf. Serv. Oper. Logist. Inform., Beijing: IEEE; 2008, p. 1806–11. https://doi.org/10.1109/SOLI.2008.4682823.

[25] Barnhart C, Hane CA, Vance PH. Using Branch-and-Price-and-Cut to Solve Origin-Destination Integer Multicommodity Flow Problems. Oper Res 2000;48:318–26. https://doi.org/10.1287/opre.48.2.318.12378.

[26] Retvdri G, Biro JJ, Cinkler T. A novel Lagrangian-relaxation to the minimum cost multicommodity flow problem and its application to OSPF traffic engineering. Proc. ISCC 2004 Ninth Int. Symp. Comput. Commun. IEEE Cat No04TH8769, Alexandria, Egypt: IEEE; 2004, p. 957–62. https://doi.org/10.1109/ISCC.2004.1358664.

[27] Moshref-Javadi M, Lee S. The customer-centric, multi-commodity vehicle routing problem with split delivery. Expert Syst Appl 2016;56:335–48. https://doi.org/10.1016/j.eswa.2016.03.030.

[28] Preskill J. Quantum Computing in the NISQ era and beyond. Quantum 2018;2:79. https://doi.org/10.22331/q-2018-08-06-79.

[29] Torlai G, Melko RG. Machine-Learning Quantum States in the NISQ Era. Annu Rev Condens Matter Phys 2020;11:325–44. https://doi.org/10.1146/annurev-conmatphys-031119-050651.

[30] Progress reports: interpreting the node log - IBM Documentation 2021. https://www.ibm.com/docs/en/icos/20.1.0?topic=mip-progress-reports-interpreting-node-log (accessed January 24, 2024).